\documentclass[12pt]{article}
\usepackage{amssymb,amsmath,amsthm,amsfonts,enumerate, bbm, mathdots}

\usepackage{latexsym}
\usepackage{amscd}
\usepackage{stmaryrd}
\usepackage{color}
\usepackage[all]{xypic}
\usepackage{epsfig}
\usepackage{graphics}
\usepackage{ifthen}
\usepackage{varioref}
\usepackage{rotating}
\usepackage{extarrows}
\usepackage{cite}
\usepackage{mathrsfs}

\numberwithin{equation}{section}

\theoremstyle{plain}
\newtheorem{thm}{Theorem}[section]
\newtheorem{cor}[thm]{Corollary}
\newtheorem{lem}[thm]{Lemma}
\newtheorem{prop}[thm]{Proposition}
\newtheorem{rem}[thm]{Remark}

\newmuskip\pFqmuskip
\newcommand*\pFq[6][8]{%
  \begingroup 
  \pFqmuskip=#1mu\relax
  \mathchardef\normalcomma=\mathcode`,
  \mathcode`\,=\string"8000
  \begingroup\lccode`\~=`\,
  \lowercase{\endgroup\let~}\pFqcomma
  {}_{#2}F_{#3}{\left(\genfrac..{0pt}{}{#4}{#5};#6\right)}%
  \endgroup
}
\newcommand{\pFqcomma}{{\normalcomma}\mskip\pFqmuskip}



\definecolor{darkgreen}{rgb}{0.0625,0.64,0.0625}
\usepackage[
            pdfstartview=FitH,
            CJKbookmarks=true,
            bookmarksnumbered=true,
            bookmarksopen=true,
            colorlinks,
            pdfborder=001,
            linkcolor=blue,
            anchorcolor=green,
            citecolor=red
            ]{hyperref}

\def\proof{\noindent {\bf Proof.\;}}

\def\dep{\operatorname{dep}}

\allowdisplaybreaks

\begin{document}

\title{Some evaluations of interpolated multiple zeta values and interpolated multiple $t$-values}
\date{~}

\date{\small ~ \qquad\qquad School of Mathematical Sciences\newline Key Laboratory of Intelligent Computing and Applications (Ministry of Education)\newline Tongji University, Shanghai 200092, China}

\author{Zhonghua Li\thanks{E-mail address: zhonghua\_li@tongji.edu.cn} ~and ~Zhenlu Wang\thanks{E-mail address: zlw@tongji.edu.cn}}

\maketitle

\begin{abstract}
In this paper, we study the evaluation formulas of the interpolated multiple zeta values and the interpolated multiple $t$-values with indices involving $1,2,3$. To get these evaluations, we derive the corresponding algebraic relations in the harmonic algebra.
\end{abstract}

{\small

}

\section{Introduction}
A finite sequence of positive integers is called an index. An index $\mathbf{k}=(k_1,\ldots,k_n)$ is admissible if $k_1\geq2$. The depth of $\mathbf{k}$ is defined by $\dep(\mathbf{k})=n$.

To study the multiple zeta values and the multiple zeta-star values simultaneously, S. Yamamoto introuced the interpolated multiple zeta values in \cite{Yamamoto}, which are interpolation polynomials of these two types of values. Let $r$ be a variable. For an admissible index $\mathbf{k}=(k_1,\ldots,k_n)$, the interpolated multiple zeta value ($r$-MZV) $\zeta^r(\mathbf{k})$ is defined by
\begin{align*}
    \zeta^r(\mathbf{k})=\zeta^r(k_1,\ldots,k_n)=\sum\limits_{\mathbf{p}=(k_1\Box k_2\Box\cdots\Box k_n)\atop \Box=``,"\text{or}``+"}r^{n-\dep(\mathbf{p})}\zeta(\mathbf{p}),
\end{align*}
where for an admissible index $\mathbf{k}=(k_1,\ldots,k_n)$, the multiple zeta value (MZV) $\zeta(\mathbf{k})$ is defined by
$$\zeta(\mathbf{k})=\zeta(k_1,\ldots,k_n)=\sum\limits_{m_1>\cdots>m_n>0}\frac{1}{m_1^{k_1}\cdots m_n^{k_n}}.$$
It is obvious that $\zeta^0(\mathbf{k})=\zeta(\mathbf{k})$. And one can find that
$$\zeta^1(\mathbf{k})=\zeta^\star(\mathbf{k})=\zeta^\star(k_1,\ldots,k_n)=\sum\limits_{m_1\geq \cdots\geq m_n>0}\frac{1}{m_1^{k_1}\cdots m_n^{k_n}},$$
which is called the multiple zeta-star value (MZSV). See \cite{Ihara-Kajikawa-Ohno-Okuda,Yamamoto,Zagier2012,Li2013} for more details.

In recent years, an odd variant of multiple zeta(-star) values, known as the multiple $t$-(star) values, has been widely studied, see \cite{Charlton,Li-Wang,Murakami2021}. For an admissible index $\mathbf{k}=(k_1,\ldots,k_n)$, the multiple $t$-value (M$t$V) and the multiple $t$-star value (M$t$SV) are defined respectively as
$$t(\mathbf{k})=t(k_1,\ldots,k_n)=\sum\limits_{m_1>\cdots>m_n>0\atop m_i:\text{odd}}\frac{1}{m_1^{k_1}\cdots m_n^{k_n}}$$
and
$$t^\star(\mathbf{k})=t(k_1,\ldots,k_n)=\sum\limits_{m_1\geq\cdots\geq m_n>0\atop m_i:\text{odd}}\frac{1}{m_1^{k_1}\cdots m_n^{k_n}}.$$
Similarly, for an admissible index $\mathbf{k}=(k_1,\ldots,k_n)$, the interpolated multiple $t$-value ($r$-M$t$V) is defined in \cite{Hoffman-Charlton} by
\begin{align*}
t^r(\mathbf{k})=t^r(k_1,\ldots,k_n)=\sum\limits_{\mathbf{p}=(k_1\Box k_2\Box\cdots\Box k_n)\atop \Box=``,"\text{or}``+"}r^{n-\dep(\mathbf{p})}t(\mathbf{p}).
\end{align*}
It is easy to see that $t^0(\mathbf{k})=t(\mathbf{k})$ and $t^1(\mathbf{k})=t^\star(\mathbf{k})$.

Multiple zeta(-star) values and multiple $t$-(star) values with indices involving $1,2,3$ are widely studied. D. Zagier gave evaluation formulas for $\zeta(2,\ldots,2,3,2,\ldots,2)$ and $\zeta^\star(2,\ldots,2,3,2,\ldots,2)$ in \cite{Zagier2012}. The analogue evaluation formulas for $t(2,\ldots,2,3,2,\ldots,2)$ and $t^\star(2,\ldots,2,3,2,\ldots,2)$ were showed in \cite{Murakami2021} and \cite{Li-Wang} respectively. S. Charlton gave an evaluation formula for regularized multiple $t$-values $t^V_{\ast}(2,\ldots,2,1,2,\ldots,2)$ in \cite{Charlton} and the corresponding formula for regularized multiple $t$-star values were showed in \cite{Li-Wang}. 

In this paper, we discuss evaluation formulas of $r$-MZVs and $r$-M$t$Vs. In Section \ref{Sec:algebra}, we use the harmonic algebra to derive some algebraic relations. As applications, we show some evaluation formulas of $r$-MZVs and $r$-M$t$Vs with indices involving $1,2,3$ in Section \ref{Sec:applications}.

\section{Harmonic algebra}\label{Sec:algebra}
\subsection{Algebraic setup}
Here we recall the harmonic algebra introduced in \cite{Ihara-Kajikawa-Ohno-Okuda,Yamamoto}.

Let $\mathfrak{A}$ be a commutative $\mathbb{Q}$-algebra and $A$ be an alphabet with non-commutative letters. Denote by $\mathfrak{h}^1$ the non-commutative polynomial algebra generated by the set $A$ over the algebra $\mathfrak{A}$, and $\mathfrak{z}$ the $\mathfrak{A}$-submodule of $\mathfrak{h}^1$ generated by $A$.

Assume the $\mathfrak{z}$ admits a commutative $\mathfrak{A}$-algebra structure (not necessary unitary) given by a product $\circ$, which is called the circle product. Define an action of $\mathfrak{z}$ on $\mathfrak{h}^1$ by $\mathfrak{A}$-linearity and the rules:
\begin{align*}
  a\circ1_w=0\quad\text{and}\quad a\circ(bw)=(a\circ b)w,
\end{align*}
where $1_w$ is the empty word, $a,b\in A$ and $w\in\mathfrak{h}^1$ is a word. We denote $\overbrace{a\circ a\circ\cdots\circ a}^n$ by $a^{\circ n}$ for convenience.

The harmonic product $\ast$ on $\mathfrak{h}^1$ is defined by $\mathfrak{A}$-bilinearity and the rules:
\begin{align*}
  &1_w\ast w=w\ast 1_w=w,\\
  &aw_1\ast bw_2=a(w_1\ast bw_2)+b(aw_1\ast w_2)+(a\circ b)(w_1\ast w_2),
\end{align*}
for all $a,b\in A$ and all words $w,w_1,w_2\in\mathfrak{h}^1$. Under the harmonic product, $\mathfrak{h}^1$ becomes a unitary commutative $\mathfrak{A}$-algebra which is called the harmonic algebra and is denoted by $\mathfrak{h}^1_\ast$.

Define an $\mathfrak{A}[r]$-linear map $\mathrm{S}^r$ on $\mathfrak{h}^1[r]$ by $\mathrm{S}^r(1_w)=1$ and
\begin{align*}
  \mathrm{\mathrm{S}}^r(aw)=a\mathrm{S}^r(w)+ra\circ \mathrm{S}^r(w)
\end{align*}
for all $a\in A$ and any word $w\in\mathfrak{h}^1$.

\begin{rem}
Obviously, $\mathrm{S}^0$ is the identity map on $\mathfrak{h}^1$, and $\mathrm{S}^1=\mathrm{S}$ is the map in \cite[Definition 1]{Ihara-Kajikawa-Ohno-Okuda}. Also it is easy to see that
\begin{align*}
\mathrm{S}^r(a\circ w)=a\circ \mathrm{S}^r(w)
\end{align*}
for all $a\in\mathfrak{z}$ and $w\in\mathfrak{h}^1[r]$.
\end{rem}

In the next subsection, we will deduce some formulas in the harmonic algebra $\mathfrak{h}^1_\ast$. For applications, we recall the special cases concerned with multiple zeta values and multiple $t$-values. Let $\mathfrak{A}=\mathbb{Q}$, $A=\{z_k\}_{k=1}^\infty$ and
$$\mathfrak{h}^0=\mathbb{Q}\oplus\bigoplus\limits_{k=2}^\infty z_k\mathfrak{h}^1\subset\mathfrak{h}^1.$$ The circle product on $\mathfrak{z}$ is defined by
\begin{align}\label{Def:circle product}
z_k\circ z_l=z_{k+l},\quad \text{for any}\quad k,l\geq1.
\end{align}
Then $\mathfrak{h}^0$ is a subalgebra of $\mathfrak{h}^1_\ast$, which is denoted by $\mathfrak{h}^0_\ast$.

There are homomorphisms of $\mathbb{Q}$-algebras
\begin{align*}
&\mathrm{\zeta}: \mathfrak{h}_\ast^0\to\mathbb{R},\quad\mathrm{t}: \mathfrak{h}_\ast^0\to\mathbb{R},
\end{align*}
which are determined by  $\zeta(1_w)=1=\mathrm{t}(1_w)$, and 
\begin{align*}
  \zeta(z_{k_1}\cdots z_{k_n})=\zeta(k_1,\ldots,k_n),\quad \mathrm{t}(z_{k_1}\cdots z_{k_n})=t(k_1,\ldots,k_n),
\end{align*}
where $z_{k_1}\cdots z_{k_n}\in\mathfrak{h}^0$ is a word. Similarly, we have homomorphisms of $\mathbb{Q}[r]$-algebras
\begin{align*}
\zeta^r: \mathfrak{h}_\ast^0\to\mathbb{R}[r],\quad \mathrm{t}^r: \mathfrak{h}_\ast^0\to\mathbb{R}[r],
\end{align*}
such that $\zeta^r(1_w)=1=\mathrm{t}^r(1_w)$, and for any $z_{k_1}\cdots z_{k_n}\in\mathfrak{h}^0$,
\begin{align*}
  \zeta^r(z_{k_1}\cdots z_{k_n})=\zeta^r(k_1,\ldots,k_n),\quad \mathrm{t}^r(z_{k_1}\cdots z_{k_n})=t^r(k_1,\ldots,k_n).
  \end{align*}
Then one can easily find that for any $w\in\mathfrak{h}^0$, 
\begin{align*}
\zeta^r(w)=\zeta(\mathrm{S}^r(w))
\end{align*}
and
\begin{align*}
  \mathrm{t}^r(w)=\mathrm{t}(\mathrm{S}^r(w)).
\end{align*}

The maps $\zeta$ and $\mathrm{t}$ can be extended to be homomorphisms of $\mathbb{Q}$-algebras
\begin{align*}
  \mathrm{\zeta}^T_\ast:\mathfrak{h}_{\ast}^1\to\mathbb{R}[T],\quad \mathrm{t}^V_\ast:\mathfrak{h}_{\ast}^1\to\mathbb{R}[V],
\end{align*}
such that $\mathrm{\zeta}^T_\ast(1_w)=1=\mathrm{t}^V_\ast(1_w)$, $\mathrm{\zeta}^T_\ast(z_1)=T,\mathrm{t}^V_\ast(z_1)=V$, and for any $z_{k_1}\cdots z_{k_n}\in\mathfrak{h}^0$,
\begin{align*}
  \mathrm{\zeta}^T_\ast(z_{k_1}\cdots z_{k_n})=\zeta(k_1,\ldots,k_n),\quad\mathrm{t}^V_\ast(z_{k_1}\cdots z_{k_n})=t(k_1,\ldots,k_n).
\end{align*}
Then for any index $(k_1,\ldots,k_n)$, one can define 
\begin{align*}
  \zeta^T_{\ast}(k_1,\ldots,k_n)=\mathrm{\zeta}^T_\ast(z_{k_1}\cdots z_{k_n}), \quad t^V_{\ast}(k_1,\ldots,k_n)=\mathrm{t}^V_\ast(z_{k_1}\cdots z_{k_n}),
\end{align*}
which are called the harmonic regularized multiple zeta values and the harmonic regularized multiple $t$-values respectively. Hence we have $\zeta^T_\ast(1)=\mathrm{\zeta}_\ast^T(z_1)=T$ and $t^V_\ast(1)=\mathrm{t}_\ast^V(z_1)=V$. For more details, see \cite{Charlton,Li2022,Ihara-Kaneko-Zagier}.

Furthermore, there are homomorphisms of $\mathbb{Q}[r]$-algebras
$\mathrm{\zeta}^{T,r}_\ast:\mathfrak{h}^1_\ast\to\mathbb{R}[T,r]$ and $\mathrm{t}^{V,r}_\ast:\mathfrak{h}^1_\ast\to\mathbb{R}[V,r]$, such that for any $w\in\mathfrak{h}^1$,
\begin{align*}
  \mathrm{\zeta}^{T,r}_\ast(w)=\mathrm{\zeta}^T_\ast(\mathrm{S}^r(w)),\quad \mathrm{t}^{V,r}_\ast(w)=\mathrm{t}^V_\ast(\mathrm{S}^r(w)).
\end{align*}
Similarly, for any index $(k_1,\ldots,k_n)$, define 
\begin{align*}
  \zeta^{T,r}_{\ast}(k_1,\ldots,k_n)=\mathrm{\zeta}^{T,r}_\ast(z_{k_1}\cdots z_{k_n}),\quad t^{V,r}_{\ast}(k_1,\ldots,k_n)=\mathrm{t}^{V,r}_\ast(z_{k_1}\cdots z_{k_n}),
\end{align*}
which are called the harmonic regularized interpolated multiple zeta values and the harmonic regularized interpolated multiple $t$-values respectively. 

Moreover, we denote $\zeta^{T,1}_{\ast}(k_1,\ldots,k_n)$ by $\zeta^{T,\star}_{\ast}(k_1,\ldots,k_n)$, and call them the harmonic regularized multiple zeta-star values. In a similar way, we call $t^{V,1}_{\ast}(k_1,\ldots,k_n)=t^{V,\star}_{\ast}(k_1,\ldots,k_n)$ the harmonic regularized multiple $t$-star values.

\subsection{Algebraic relations}
We give the following lemma which generalizes \cite[Lemma 4.5]{Li2022}.
\begin{lem}\label{Lem:S^r(a^n)}
  For any $a_1,a_2,\ldots,a_n\in\mathfrak{z}$, we have
  \begin{align}\label{Eq:recurrent relation for S^r}
    \mathrm{S}^r(a_1a_2\cdots a_n)=\sum\limits_{i=1}^nr^{i-1}(a_1\circ a_2\circ\cdots\circ a_i)\mathrm{S}^r(a_{i+1}\cdots a_n).
  \end{align}
  In particular, for any positive integer $n$ and any $a\in\mathfrak{z}$, we have
  \begin{align}\label{Eq:recurrent relation for S^r(a^n)}
  \mathrm{S}^r(a^n)=\sum\limits_{i=1}^nr^{i-1}a^{\circ i}\mathrm{S}^r(a^{n-i}).
  \end{align}
\end{lem}
\proof
We prove \eqref{Eq:recurrent relation for S^r} by induction on $n$. The case of $n=1$ is trivial. Assume that $n\geq2$, since
\begin{align*}
  \mathrm{S}^r(a_1a_2\cdots a_n)&=a_1\mathrm{S}^r(a_2\cdots a_n)+ra_1\circ \mathrm{S}^r(a_2\cdots a_n)\\
  &=a_1\mathrm{S}^r(a_2\cdots a_n)+ra_1\circ \sum\limits_{i=2}^{n}r^{i-2}(a_2\circ\cdots\circ a_i)\mathrm{S}^r(a_{i+1}\cdots a_n)\\
  &=a_1\mathrm{S}^r(a_2\cdots a_n)+\sum\limits_{i=2}^{n}r^{i-1}(a_1\circ a_2\circ\cdots\circ a_i)\mathrm{S}^r(a_{i+1}\cdots a_n),
\end{align*}
we finish the proof.
\qed

The following proposition generalizes \cite[Corollary 5.5]{Hoffman-Ihara}.
\begin{prop}\label{Prop:S^r(a^n)}
  Let $u$ be a formal variable. For any $a\in\mathfrak{z}$, we have
  \begin{align}\label{Eq:S^r-S}
    \mathrm{S}^r\left(\frac{1}{1-au}\right)=\frac{1}{1-a(1-r)u}\ast \mathrm{S}\left(\frac{1}{1-aru}\right),
  \end{align}
  or equivalently, for any positive integer $n$, we have
  \begin{align}\label{Eq:S^r(a^n)}
   \mathrm{S}^r(a^n)=\sum\limits_{i=0}^nr^{n-i}(1-r)^ia^i\ast \mathrm{S}(a^{n-i}).
  \end{align}
\end{prop}
\proof
We prove \eqref{Eq:S^r(a^n)} by induction on $n$. It is obvious that \eqref{Eq:S^r(a^n)} holds for $n=1$. Then assume that $n\geq 2$. Using \eqref{Eq:recurrent relation for S^r(a^n)} with $r=1$, the right hand side of \eqref{Eq:S^r(a^n)} is
\begin{align*}
&(1-r)^na^n+r^n\mathrm{S}(a^n)+\sum\limits_{i=1}^{n-1}\sum\limits_{j=1}^{i}r^i(1-r)^{n-i}a^{n-i}\ast a^{\circ j}\mathrm{S}(a^{i-j}),
\end{align*}
which is $r^n\mathrm{S}(a^n)+R_1+R_2+R_3$, where
\begin{align*}
R_1&=(1-r)^na^n+a\sum\limits_{i=1}^{n-1}\sum\limits_{j=1}^{i}r^i(1-r)^{n-i}a^{n-i-1}\ast a^{\circ j}\mathrm{S}(a^{i-j}),\\
R_2&=\sum\limits_{i=1}^{n-1}\sum\limits_{j=1}^{i}r^i(1-r)^{n-i}a^{\circ j}(a^{n-i}\ast \mathrm{S}(a^{i-j})),\\
R_3&=\sum\limits_{i=1}^{n-1}\sum\limits_{j=1}^{i}r^i(1-r)^{n-i}a^{\circ (j+1)}(a^{n-i-1}\ast \mathrm{S}(a^{i-j})).
\end{align*}
For $R_1$, we have 
\begin{align*}
R_1=(1-r)^na^n+a\sum\limits_{i=1}^{n-1}r^{i}(1-r)^{n-i}a^{n-i-1}\ast \mathrm{S}(a^{i}),
\end{align*}
which is $(1-r)a\mathrm{S}^r(a^{n-1})$ by the inductive hypothesis.
For $R_2$, by changing the order of summation, we have
\begin{align*}
R_2&=\sum\limits_{j=1}^{n-1}r^ja^{\circ j}\sum\limits_{i=j}^{n-1}r^{i-j}(1-r)^{n-i}a^{n-i}\ast \mathrm{S}(a^{i-j})\\
&=\sum\limits_{j=1}^{n-1}r^ja^{\circ j}\sum\limits_{i=0}^{n-j-1}r^{i}(1-r)^{n-i-j}a^{n-i-j}\ast \mathrm{S}(a^{i}).
\end{align*}
Using the inductive hypothesis, we have
\begin{align*}
R_2&=\sum\limits_{j=1}^{n-1}r^ja^{\circ j}(\mathrm{S}^r(a^{n-j})-r^{n-j}\mathrm{S}(a^{n-j}))\\
&=\sum\limits_{j=1}^{n-1}r^ja^{\circ j}\mathrm{S}^r(a^{n-j})-\sum\limits_{j=1}^{n-1}r^na^{\circ j}\mathrm{S}(a^{n-j}).
\end{align*}
Similarly, for $R_3$, we get
\begin{align*}
R_3&=\sum\limits_{j=1}^{n-1}r^ja^{\circ (j+1)}\sum\limits_{i=j}^{n-1}r^{i-j}(1-r)^{n-i}a^{n-i-1}\ast \mathrm{S}(a^{i-j})\\
&=\sum\limits_{j=1}^{n-1}(1-r)r^ja^{\circ (j+1)}\mathrm{S}^r(a^{n-j-1})\\
&=\sum\limits_{j=2}^{n}r^{j-1}a^{\circ j}\mathrm{S}^r(a^{n-j})-\sum\limits_{j=2}^{n}r^{j}a^{\circ j}\mathrm{S}^r(a^{n-j}).
\end{align*}
Therefore, the right hand side of \eqref{Eq:S^r(a^n)} is
\begin{align*}
\sum\limits_{j=1}^nr^{j-1}a^{\circ j}\mathrm{S}^r(a^{n-j}),
\end{align*}
which is $\mathrm{S}^r(a^n)$ by \eqref{Eq:recurrent relation for S^r(a^n)}. 
\qed

We may regard \eqref{Eq:S^r-S} as a generalization of
\begin{align}\label{Eq:S_1}
  \mathrm{S}\left(\frac{1}{1-au}\right)\ast\frac{1}{1+au}=1,
\end{align}
which was proved in \cite[Corollary 1]{Ihara-Kajikawa-Ohno-Okuda}. For a positive integer $k\geq 2$, let $a=z_k$ and the circle product follows \eqref{Def:circle product}. By applying the evaluation maps $\zeta$ and $\mathrm{t}$ to \eqref{Eq:S_1}, we get
\begin{align}\label{Eq:zeta and t k}
\sum\limits_{i=0}^{n}(-1)^i\zeta(\{k\}^i)\zeta^\star(\{k\}^{n-i})=\sum\limits_{i=0}^{n}(-1)^it(\{k\}^i)t^\star(\{k\}^{n-i})=\begin{cases}
1 & n=0,\\
0 & n\geq 1.
\end{cases}
\end{align}
Also, by applying the evaluation maps $\zeta$ and $\mathrm{t}$ to \eqref{Eq:S^r-S}, we have
\begin{align}
\zeta^{r}(\{k\}^n)&=\sum\limits_{i=0}^nr^{n-i}(1-r)^i\zeta(\{k\}^i)\zeta^\star(\{k\}^{n-i}),\label{Eq:zeta^r k}\\
t^{r}(\{k\}^n)&=\sum\limits_{i=0}^nr^{n-i}(1-r)^it(\{k\}^i)t^\star(\{k\}^{n-i}).\label{Eq:t^r k}
\end{align}

The following two lemmas can be deduced from \cite[Theorem 2.2]{Li2013}.
\begin{lem}
Let $u$ be a formal variable. For any $a_1,a_2\in\mathfrak{z}$, we have
  \begin{align*}
    \mathrm{S}\left(\frac{1}{1-a_1u}a_2\right)&=a_2\frac{1}{1+a_1u}\ast \mathrm{S}\left(\frac{1}{1-a_1u}\right),
  \end{align*}
  or equivalently, for any positive integer $p$, we have
  \begin{align}
  \mathrm{S}(a_1^pa_2)&=\sum\limits_{i=0}^p(-1)^ia_2a_1^i\ast \mathrm{S}(a_1^{p-i}).\label{Eq:S(a_1^pa_2)}
  \end{align}
\end{lem}

\begin{lem}
  Let $u$ be a formal variable. For any $a_1,a_2\in\mathfrak{z}$, we have
    \begin{align*}
      \mathrm{S}\left(a_1\frac{1}{1-a_2u}\right)&=\frac{1}{1+a_2u}a_1\ast \mathrm{S}\left(\frac{1}{1-a_2u}\right),
    \end{align*}
    or equivalently, for any positive integer $q$, we have
    \begin{align*}
    \mathrm{S}(a_1a_2^q)&=\sum\limits_{i=0}^q(-1)^ia_2^ia_1\ast \mathrm{S}(a_2^{q-i}).
    \end{align*}
  \end{lem}

We have the following generalized results.

\begin{prop}\label{Prop:S^r(a_1^pa_2)}
  Let $u$ be a formal variable. For any $a_1,a_2\in\mathfrak{z}$, we have
  \begin{align*}
    \mathrm{S}^r\left(\frac{1}{1-a_1u}a_2\right)=\frac{1}{1-a_1(1-r)u}a_2\frac{1}{1+a_1ru}\ast \mathrm{S}\left(\frac{1}{1-a_1ru}\right),
  \end{align*}
  or equivalently, for any positive integer $p$, we have
  \begin{align}\label{Eq:S^r(a_1^pa_2)}
  \mathrm{S}^r(a_1^pa_2)=\sum\limits_{i+j=0\atop i,j\geq0}^p(-1)^jr^{p-i}(1-r)^ia_1^ia_2a_1^j\ast \mathrm{S}(a_1^{p-i-j}).
  \end{align}
\end{prop}
\proof
We prove \eqref{Eq:S^r(a_1^pa_2)} by induction on $p$. By direct calculation, it is obvious that \eqref{Eq:S^r(a_1^pa_2)} holds when $p=1$. Assume that $p\geq 2$. The right hand side of \eqref{Eq:S^r(a_1^pa_2)} is equal to 
\begin{align*}
\sum\limits_{i=0}^p(-1)^{p-i}r^{p-i}(1-r)^ia_1^ia_2a_1^{p-i}+\sum\limits_{j=0}^{p-1}(-1)^jr^pa_2a_1^j\ast \mathrm{S}(a_1^{p-j})\\
+\sum\limits_{i=1}^{p-1}\sum\limits_{j=0}^{p-i-1}(-1)^jr^{p-i}(1-r)^ia_1^ia_2a_1^j\ast \mathrm{S}(a_1^{p-i-j}).
\end{align*}
Using \eqref{Eq:recurrent relation for S^r(a^n)} with $r=1$ and \eqref{Eq:S(a_1^pa_2)}, the right hand side of \eqref{Eq:S^r(a_1^pa_2)} becomes
\begin{align*}
\sum\limits_{i=0}^p(-1)^{p-i}r^{p-i}(1-r)^ia_1^ia_2a_1^{p-i}+r^p\mathrm{S}(a_1^pa_2)+(-1)^{p-1}r^pa_2a_1^p+Q_1+Q_2+Q_3,
\end{align*}
where 
\begin{align*}
Q_1&=a_1\sum\limits_{i=1}^{p-1}\sum\limits_{j=0}^{p-i-1}\sum\limits_{n=1}^{p-i-j}(-1)^jr^{p-i}(1-r)^ia_1^{i-1}a_2a_1^j\ast a_1^{\circ n}\mathrm{S}(a_1^{p-i-j-n}),\\
Q_2&=\sum\limits_{i=1}^{p-1}\sum\limits_{j=0}^{p-i-1}\sum\limits_{n=1}^{p-i-j}(-1)^jr^{p-i}(1-r)^ia_1^{\circ n}(a_1^ia_2a_1^j\ast \mathrm{S}(a_1^{p-i-j-n})),\\
Q_3&=\sum\limits_{i=1}^{p-1}\sum\limits_{j=0}^{p-i-1}\sum\limits_{n=1}^{p-i-j}(-1)^jr^{p-i}(1-r)^ia_1^{\circ (n+1)}(a_1^{i-1}a_2a_1^j\ast \mathrm{S}(a_1^{p-i-j-n})).
\end{align*}
For $Q_1$, we have
\begin{align*}
Q_1&=a_1\sum\limits_{i=1}^{p-1}\sum\limits_{j=0}^{p-i-1}(-1)^jr^{p-i}(1-r)^ia_1^{i-1}a_2a_1^j\ast \mathrm{S}(a_1^{p-i-j})\\
&=a_1\sum\limits_{i+j=0\atop i,j\geq0}^{p-2}(-1)^jr^{p-i-1}(1-r)^{i+1}a_1^ia_2a_1^j\ast \mathrm{S}(a_1^{p-i-j-1}).
\end{align*}
According to the inductive hypothesis, we get 
\begin{align*}
Q_1&=(1-r)a_1\mathrm{S}^r(a_1^{p-1}a_2)-\sum\limits_{i=1}^p(-1)^{p-i}r^{p-i}(1-r)^ia_1^ia_2a_1^{p-i}.
\end{align*}
For $Q_2$, changing the order of summation, we have
\begin{align*}
Q_2
&=\sum\limits_{n=1}^{p-1}r^na_1^{\circ n}\sum\limits_{i=1}^{p-n}\sum\limits_{j=0}^{p-n-i}(-1)^jr^{p-n-i}(1-r)^ia_1^ia_2a_1^j\ast \mathrm{S}(a_1^{p-i-j-n}).
\end{align*}
Using the inductive hypothesis and \eqref{Eq:S(a_1^pa_2)}, we get
\begin{align*}
Q_2&=\sum\limits_{n=1}^{p-1}r^na_1^{\circ n}\left(\mathrm{S}^r(a_1^{p-n}a_2)-\sum\limits_{j=0}^{p-n}(-1)^jr^{p-n}a_2a_1^j\ast \mathrm{S}(a_1^{p-n-j})\right)\\
&=\sum\limits_{n=1}^{p-1}r^na_1^{\circ n}\mathrm{S}^r(a_1^{p-n}a_2)-\sum\limits_{n=1}^{p-1}r^pa_1^{\circ n}\sum\limits_{j=0}^{p-n}(-1)^ja_2a_1^j\ast \mathrm{S}(a_1^{p-n-j})\\
&=\sum\limits_{n=1}^{p-1}r^na_1^{\circ n}\mathrm{S}^r(a_1^{p-n}a_2)-\sum\limits_{n=1}^{p-1}r^pa_1^{\circ n}\mathrm{S}(a_1^{p-n}a_2).
\end{align*}
Similarly, for $Q_3$, we have
\begin{align*}
Q_3&=\sum\limits_{n=1}^{p-1}r^na_1^{\circ (n+1)}\sum\limits_{i=1}^{p-n}\sum\limits_{j=0}^{p-n-i}(-1)^jr^{p-n-i}(1-r)^ia_1^{i-1}a_2a_1^j\ast \mathrm{S}(a_1^{p-i-j-n})\\
&=\sum\limits_{n=1}^{p-1}(1-r)r^na_1^{\circ (n+1)}\mathrm{S}^r(a_1^{p-n-1}a_2)\\
&=\sum\limits_{n=2}^{p}r^{n-1}a_1^{\circ n}\mathrm{S}^r(a_1^{p-n}a_2)-\sum\limits_{n=2}^{p}r^{n}a_1^{\circ n}\mathrm{S}^r(a_1^{p-n}a_2).
\end{align*}
Hence, the right hand side of \eqref{Eq:S^r(a_1^pa_2)} is
\begin{align*}
r^p\mathrm{S}(a_1^pa_2)-r^p\sum\limits_{n=1}^pa_1^{\circ n}\mathrm{S}(a_1^{p-n}a_2)+\sum\limits_{n=1}^pr^{n-1}a_1^{\circ n}\mathrm{S}^r(a_1^{p-n}a_2),
\end{align*}
which is $\mathrm{S}^r(a_1^pa_2)$ by \eqref{Eq:recurrent relation for S^r}.
\qed

\begin{prop}\label{Prop:S^r(a_1a_2^q)}
  Let $u$ be a formal variable. For any $a_1,a_2\in\mathfrak{z}$, we have
  \begin{align*}
    \mathrm{S}^r\left(a_1\frac{1}{1-a_2u}\right)=\frac{1}{1+a_2ru}a_1\frac{1}{1-a_2(1-r)u}\ast \mathrm{S}\left(\frac{1}{1-a_2ru}\right),
  \end{align*}
  or equivalently, for any positive integer $q$, we have
  \begin{align}\label{Eq:S^r-r}
  \mathrm{S}^r(a_1a_2^q)=\sum\limits_{i+j=0\atop i,j\geq0}^q(-1)^ir^{q-j}(1-r)^ja_2^ia_1a_2^j\ast \mathrm{S}(a_2^{q-i-j}).
  \end{align}
\end{prop}
\proof
We prove \eqref{Eq:S^r-r} by induction on $q$. The case of $q=1$ holds by direct calculation. Assume that $q\geq2$. By using \eqref{Eq:recurrent relation for S^r(a^n)} with $r=1$, the right hand side of \eqref{Eq:S^r-r} is
\begin{align*}
\sum\limits_{i=0}^q(-1)^ir^i(1-r)^{q-i}a_2^ia_1a_2^{q-i}+\sum\limits_{j=0}^{q-1}r^{q-j}(1-r)^ja_1a_2^j\ast \mathrm{S}(a_2^{q-j})+T_1+T_2+T_3,
\end{align*}
where
\begin{align*}
T_1&=a_2\sum\limits_{i=1}^{q-1}\sum\limits_{j=0}^{q-i-1}\sum\limits_{n=1}^{q-i-j}(-1)^ir^{q-j}(1-r)^ja_2^{i-1}a_1a_2^j\ast a_2^{\circ n}\mathrm{S}(a_2^{q-i-j-n}),\\
T_2&=\sum\limits_{i=1}^{q-1}\sum\limits_{j=0}^{q-i-1}\sum\limits_{n=1}^{q-i-j}(-1)^ir^{q-j}(1-r)^ja_2^{\circ n}(a_2^{i}a_1a_2^j\ast \mathrm{S}(a_2^{q-i-j-n})),\\
T_3&=\sum\limits_{i=1}^{q-1}\sum\limits_{j=0}^{q-i-1}\sum\limits_{n=1}^{q-i-j}(-1)^ir^{q-j}(1-r)^ja_2^{\circ (n+1)}(a_2^{i-1}a_1a_2^j\ast \mathrm{S}(a_2^{q-i-j-n})).
\end{align*}
For $T_1$, we have
\begin{align*}
T_1&=a_2\sum\limits_{i=1}^{q-1}\sum\limits_{j=0}^{q-i-1}(-1)^ir^{q-j}(1-r)^ja_2^{i-1}a_1a_2^j\ast \mathrm{S}(a_2^{q-i-j})\\
&=-a_2\sum\limits_{i+j=0\atop i,j\geq0}^{q-2}(-1)^ir^{q-j}(1-r)^ja_2^{i}a_1a_2^j\ast \mathrm{S}(a_2^{q-i-j-1}).
\end{align*}
Then by the inductive hypothesis, we get
\begin{align*}
T_1&=-ra_2\mathrm{S}^r(a_1a_2^{q-1})-\sum\limits_{i=1}^q(-1)^ir^i(1-r)^{q-i}a_2^ia_1a_2^{q-i}.
\end{align*}
For $T_2$, changing the order of the summation, we have
\begin{align*}
T_2&=\sum\limits_{n=1}^{q-1}r^na_2^{\circ n}\sum\limits_{i=1}^{q-n}\sum\limits_{j=0}^{q-i-n}(-1)^ir^{q-n-j}(1-r)^ja_2^ia_1a_2^j\ast \mathrm{S}(a_2^{q-n-i-j})\\
&=\sum\limits_{n=1}^{q-1}r^na_2^{\circ n}\mathrm{S}^r(a_1a_2^{q-n})-\sum\limits_{n=1}^{q-1}r^na_2^{\circ n}\sum\limits_{j=0}^{q-n}r^{q-n-j}(1-r)^ja_1a_2^j\ast \mathrm{S}(a_2^{q-n-j}).
\end{align*}
Similarly, for $T_3$, we get
\begin{align*}
T_3&=\sum\limits_{n=1}^{q-1}r^{n+1}a_2^{\circ (n+1)}\sum\limits_{i=1}^{q-n}\sum\limits_{j=0}^{q-n-i}(-1)^ir^{q-n-j-1}(1-r)^ja_2^{i-1}a_1a_2^j\ast \mathrm{S}(a_2^{q-n-i-j})\\
&=\sum\limits_{n=1}^{q-1}r^{n+1}a_2^{\circ (n+1)}\sum\limits_{i=0}^{q-n-1}\sum\limits_{j=0}^{q-n-i-1}(-1)^{i-1}r^{q-n-j-1}(1-r)^ja_2^{i}a_1a_2^j\ast \mathrm{S}(a_2^{q-n-i-j-1})\\
&=-\sum\limits_{n=2}^qr^na_2^{\circ n}\mathrm{S}^r(a_1a_2^{q-n}).
\end{align*}
Therefore, the right hand side of \eqref{Eq:S^r-r} becomes
\begin{align*}
&\sum\limits_{j=0}^qr^{q-j}(1-r)^ja_1a_2^j\ast \mathrm{S}(a_2^{q-j})-\sum\limits_{n=1}^q\sum\limits_{j=0}^{q-n}r^{q-j}(1-r)^ja_2^{\circ n}(a_1a_2^j\ast \mathrm{S}(a_2^{q-n-j}))\\
&=\sum\limits_{j=0}^q\sum\limits_{n=1}^{q-j}r^{q-j}(1-r)^j\left\{a_1(a_2^j\ast a_2^{\circ n}\mathrm{S}(a_2^{q-n-j}))+(a_1\circ a_2^{\circ n})(a_2^j\ast \mathrm{S}(a_2^{q-n-j})\right\}.
\end{align*}
Then by \eqref{Eq:recurrent relation for S^r(a^n)}, \eqref{Eq:S^r(a^n)} and the inductive hypothesis, we find that the right hand side of \eqref{Eq:S^r-r} is equal to
\begin{align*}
a_1\mathrm{S}^r({a_2^q})+\sum\limits_{n=1}^qr^n(a_1\circ a_2^{\circ n})\mathrm{S}(a_2^{q-n}),
\end{align*}
which is $\mathrm{S}^r(a_1a_2^q)$.
\qed

\section{Evaluation formulas}\label{Sec:applications}

In the special cases with multiple zeta values and multiple $t$-values, we deduce some evaluation formulas from Proposition \ref{Prop:S^r(a_1^pa_2)} and Proposition \ref{Prop:S^r(a_1a_2^q)}. 

 For any positive integer $n$, we give the following notations
\begin{align*}
\Phi_r^n&=\frac{2(r-1)^n-(\sqrt{r-1}+\sqrt{r})^{2n}-(\sqrt{r-1}-\sqrt{r})^{2n}}{2r},\\
\Psi_r^n&=\frac{2r^n-(\sqrt{r-1}+\sqrt{r})^{2n}-(\sqrt{r-1}-\sqrt{r})^{2n}}{2(r-1)},\\
\Theta_r^n&=\frac{(\sqrt{r-1}+\sqrt{r})^{2n}-(\sqrt{r-1}-\sqrt{r})^{2n}}{2\sqrt{r(r-1)}},\\
\Omega_r^{n}&=\frac{(\sqrt{r-1}+\sqrt{r})^{2n}+(\sqrt{r-1}-\sqrt{r})^{2n}}{2}.
\end{align*}

Considering the interpolated multiple zeta values with indices involving 2 and 3, we get the following theorem.
\begin{thm}\label{Thm:r-mzvs2-3}
For any nonnegative integers $p$ and $q$, we have
\begin{align}
\zeta^r(\{2\}^p,3)&=2\sum\limits_{n=0}^p\Phi_r^{n+1}\zeta^r(\{2\}^{p-n})\zeta(2n+3)+2(1-2^{-2p-2})\Theta_r^{p+1}\zeta(2p+3),\label{Eq:zeta(2-3)}\\
\zeta^r(3,\{2\}^q)&=2\sum\limits_{n=0}^q(1-2^{-2n-2})\Theta_r^{n+1}\zeta^r(\{2\}^{q-n})\zeta(2n+3)+2\Psi_r^{q+1}\zeta(2q+3).\label{Eq:zeta(3-2)}
\end{align}
\end{thm}
\proof
We prove \eqref{Eq:zeta(2-3)}. It is obvious that \eqref{Eq:zeta(2-3)} holds when $p=0$. For $p\geq 1$, setting $a_1=z_2$ and $a_2=z_3$ in \eqref{Eq:S^r(a_1^pa_2)}, and  applying the map $\zeta$, we have
\begin{align*}
\zeta^r(\{2\}^p,3)=\sum\limits_{i+j=0\atop i,j\geq0}^p(-1)^jr^{p-i}(1-r)^i\zeta(\{2\}^i,3,\{2\}^j)\zeta^\star(\{2\}^{p-i-j}).
\end{align*}
The following evaluation formula was showed in \cite[Theorem 1]{Zagier2012}
\begin{align}\label{Eq:zagier evaluation}
&\zeta(\{2\}^i,3,\{2\}^j)\notag\\&=2\sum\limits_{n=1}^{i+j+1}(-1)^{n}\left\{\binom{2n}{2j+2}-(1-2^{-2n})\binom{2n}{2i+1}\right\}\zeta(\{2\}^{i+j-n+1})\zeta(2n+1).
\end{align}
Then we have
\begin{align*}
\zeta^r(\{2\}^p,3)=2\sum\limits_{n=1}^{p+1}(-1)^n\zeta(2n+1)B_1(n,p)+2\sum\limits_{n=1}^{p+1}(-1)^n(1-2^{-2n})\zeta(2n+1)B_2(n,p),
\end{align*}
where
\begin{align*}
B_1(n,p)&=\sum\limits_{i+j=n-1\atop i,j\geq0}^p(-1)^{j}r^{p-i}(1-r)^i\binom{2n}{2j+2}\zeta(\{2\}^{i+j-n+1})\zeta^\star(\{2\}^{p-i-j}),\\
B_2(n,p)&=\sum\limits_{i+j=n-1\atop i,j\geq0}^p(-1)^{j-1}r^{p-i}(1-r)^i\binom{2n}{2i+1}\zeta(\{2\}^{i+j-n+1})\zeta^\star(\{2\}^{p-i-j}).
\end{align*}
For $B_1(n,p)$, we have
\begin{align*}
&B_1(n,p)\\&=\sum\limits_{j=0}^{n-1}\sum\limits_{i=n-1-j}^{p-j}(-1)^jr^{p-i}(1-r)^i\binom{2n}{2j+2}\zeta(\{2\}^{i+j-n+1})\zeta^\star(\{2\}^{p-i-j})\\
&=\sum\limits_{j=0}^{n-1}(-1)^jr^j(1-r)^{n-j-1}\binom{2n}{2j+2}\sum\limits_{i=0}^{p-n+1}r^{p-n+1-i}(1-r)^i\zeta(\{2\}^i)\zeta^\star(\{2\}^{p-n+1-i}).
\end{align*}
Using \eqref{Eq:zeta^r k}, we have
\begin{align*}
B_1(n,p)=\sum\limits_{j=0}^{n-1}(-1)^{j}r^{j}(1-r)^{n-1-j}\binom{2n}{2j+2}\zeta^r(\{2\}^{p-n+1}).
\end{align*}
Using the combinatorial identity \cite[(1.87)]{Gould}
\begin{align}\label{Eq:Gouldbinomal2}
\sum\limits_{k=0}^{n}\binom{2n}{2k}x^k=\frac{(1+\sqrt{x})^{2n}+(1-\sqrt{x})^{2n}}{2}\quad\text{for}\quad n\geq0,
\end{align}
we get
\begin{align*}
B_1(n,p)=(-1)^n\Phi_r^n\zeta^r(\{2\}^{p-n+1}).
\end{align*}
For $B_2(n,p)$, we have
\begin{align*}
B_2(n,p)&=\sum\limits_{i=0}^{n-1}\sum\limits_{j=n-1-i}^{p-i}(-1)^{j-1}r^{p-i}(1-r)^i\binom{2n}{2i+1}\zeta(\{2\}^{i+j-n+1})\zeta^\star(\{2\}^{p-i-j})\\
&=\sum\limits_{i=0}^{n-1}(-1)^{n-i}r^{p-i}(1-r)^{i}\binom{2n}{2i+1}\sum\limits_{j=0}^{p-n+1}(-1)^j\zeta(\{2\}^j)\zeta^\star(\{2\}^{p-n+1-j}).
\end{align*}
Using \eqref{Eq:zeta and t k} and the combinatorial identity \cite[(1.95)]{Gould}
\begin{align}\label{Eq:Gouldbinomal1}
\sum\limits_{k=0}^{n-1}\binom{2n}{2k+1}x^k=\frac{(1+\sqrt{x})^{2n}-(1-\sqrt{x})^{2n}}{2\sqrt{x}}\quad\text{for}\quad n\geq1,
\end{align}
we get
\begin{align*}
B_2(n,p)=\begin{cases}
(-1)^{p-1}\Theta_r^{p+1} & n=p+1,\\
0 & 1\leq n\leq p.
\end{cases}
\end{align*}
Therefore, we have
\begin{align*}
\zeta^r(\{2\}^p,3)=2\sum\limits_{n=1}^{p+1}\zeta(2n+1)\Phi_r^n\zeta^r(\{2\}^{p-n+1})+2(1-2^{-2p-2})\Theta_r^{p+1}\zeta(2p+3),
\end{align*}
which completes the proof.
\qed

Taking $r=0$, \eqref{Eq:zeta(2-3)} and \eqref{Eq:zeta(3-2)} are special cases of \eqref{Eq:zagier evaluation} when $j=0$ and $i=0$. Taking $r=1$, \eqref{Eq:zeta(2-3)} and \eqref{Eq:zeta(3-2)} are special cases of \cite[(45)]{Zagier2012}. Furthermore, we have the following evaluation formulas when $r=\frac{1}{2}$.

\begin{cor}
For any nonnegative integers $p$ and $q$, we have
\begin{align*}
\zeta^{\frac{1}{2}}(\{2\}^p,3)&=4\sum\limits_{n=0}^p(-1)^{n-1}\left(2^{-n-1}+\sin\frac{n\pi}{2}\right)\zeta^\frac{1}{2}(\{2\}^{p-n})\zeta(2n+3)\\&\quad+4(1-2^{-2p-2})\cos\frac{p\pi}{2}\zeta(2p+3),\\
\zeta^{\frac{1}{2}}(3,\{2\}^q)&=4\sum\limits_{n=0}^q(1-2^{-2n-2})\cos\frac{n\pi}{2}\zeta^{\frac{1}{2}}(\{2\}^{q-n})\zeta(2n+3)\\&\quad-4\left(2^{-q-1}+\sin\frac{q\pi}{2}\right)\zeta(2q+3).
\end{align*}
\end{cor}
\proof 
Since $\Theta_{1/2}^n=2\sin\frac{n\pi}{2}$ and  
\begin{align}\label{Eq:zeta-two-one-two}
  \Phi_{1/2}^n=(-1)^{n-1}\Psi_{1/2}^n=2(-1)^n\left(2^{-n}-\cos\frac{n\pi}{2}\right), 
\end{align}
we get the desired results.
\qed

In the following, we consider the harmonic regularized multiple zeta values with indices involving 1 and 2. We denote $\delta_{a,b}$ the Kronecker symbol. First, we give the following lemma. 

\begin{lem}\label{Lem:rmzvs-two-one}
  For any nonnegative integers $p$ and $q$, we have
  \begin{align}\label{Eq:rmzvs-two-one-two}
    &\zeta_{\ast}^T(\{2\}^p,1,\{2\}^q)\notag\\&=2\sum\limits_{n=1}^{p+q}(-1)^n\left\{\binom{2n}{2p}-(1-2^{-2n})\binom{2n}{2q+1}\right\}\zeta(\{2\}^{p+q-n})\zeta(2n+1)+\delta_{p,0}T\zeta(\{2\}^q).
  \end{align}
\end{lem}
\proof
If $p\geq1$, by the duality formula (see \cite{Zagier1992})
\begin{align*}
  \zeta(\{2\}^p,1,\{2\}^q)=\zeta(\{2\}^q,3,\{2\}^{q-1}),
\end{align*}
we get \eqref{Eq:zeta-two-one-two} from \eqref{Eq:zagier evaluation} directly. If $p=0$, by the harmonic product, we have
\begin{align*}
  &\zeta^T_{\ast}(1,\{2\}^q)\\&=T\zeta(\{2\}^q)-\sum\limits_{i=1}^q\zeta(\{2\}^i,1,\{2\}^{q-i})-\sum\limits_{i=0}^{q-1}\zeta(\{2\}^i,3,\{2\}^{q-1-i})\\
  &=T\zeta(\{2\}^q)-2\sum\limits_{i=1}^q\zeta(\{2\}^i,1,\{2\}^{q-i})\\
  &=T\zeta(\{2\}^q)-4\sum\limits_{n=1}^q(-1)^n\sum\limits_{i=1}^q\left(\binom{2n}{2i}-(1-2^{-2n})\binom{2n}{2q-2i+1}\right)\zeta(\{2\}^{q-n})\zeta(2n+1).
\end{align*}
Since
\begin{align*}
  \sum\limits_{i=1}^q\left\{\binom{2n}{2i}-(1-2^{-2n})\binom{2n}{2q-2i+1}\right\}=2^{2n-1}-1-(1-2^{-2n})2^{2n-1}=-\frac{1}{2},
\end{align*}
we get 
\begin{align*}
  \zeta^T_{\ast}(1,\{2\}^q)&=T\zeta(\{2\}^q)+2\sum\limits_{n=1}^q(-1)^n\zeta(\{2\}^{q-n})\zeta(2n+1).
\end{align*}
Thus, we finish the proof.
\qed

Then we get evaluations for the harmonic regularized interpolated multiple zeta values.
\begin{thm}\label{Thm:r-mzvs2-1}
  For any nonnegative integers $p$ and $q$, we have
  \begin{align}
    \zeta_{\ast}^{T,r}(\{2\}^p,1)=2\sum\limits_{n=1}^{p}(1-r)(1-2^{-2n})\Theta_r^n\zeta^r(\{2\}^{p-n})\zeta(2n+1)+2\Omega_r^p\zeta^T_\ast(2p+1)-\delta_{p,0}T,\label{Eq:rmzvs-two-one}\\
    \zeta_{\ast}^{T,r}(1,\{2\}^q)=2\sum\limits_{n=1}^{q}\Omega_r^n\zeta^r(\{2\}^{q-n})\zeta(2n+1)-2r(1-2^{-2q})\Theta_r^q\zeta^T_\ast(2q+1)+T\zeta^r(\{2\}^q).\label{Eq:rmzvs-one-two}
  \end{align}
\end{thm}
\proof
We prove \eqref{Eq:rmzvs-two-one}. Since $\Omega_r^0=1$, it is easy to see that \eqref{Eq:rmzvs-two-one} holds when $p=0$. For $p\geq 1$, let $a_1=z_2$ and $a_2=z_1$ in \eqref{Eq:S^r-r}, then by applying the map $\mathrm{\zeta}^T_\ast$, we have 
\begin{align*}
  \zeta^{T,r}_\ast(\{2\}^p,1)=\sum\limits_{i+j=0\atop i,j\geq0}^p(-1)^jr^{p-i}(1-r)^i\zeta^T_\ast(\{2\}^i,1,\{2\}^j)\zeta^\star(\{2\}^{p-i-j}).
\end{align*}
Using Lemma \ref{Lem:rmzvs-two-one}, we have
\begin{align*}
  &\zeta^{T,r}_\ast(\{2\}^p,1)\\
  &=T\sum\limits_{i+j=0\atop i,j\geq0}^p(-1)^jr^{p-i}(1-r)^i\delta_{i,0}\zeta(\{2\}^j)\zeta^\star(\{2\}^{p-i-j})\\
  &\quad+2\sum\limits_{i+j=0\atop i,j\geq0}^p\sum\limits_{n=1}^{i+j}(-1)^{n+j}r^{p-i}(1-r)^i\left\{\binom{2n}{2i}-(1-2^{-2n})\binom{2n}{2j+1}\right\}\\
  &\quad\quad\qquad\qquad\times\zeta(\{2\}^{i+j-n})\zeta^\star(\{2\}^{p-i-j})\zeta(2n+1).
\end{align*}
As $p\geq1$, by \eqref{Eq:zeta and t k}, we get 
\begin{align*}
  T\sum\limits_{i+j=0\atop i,j\geq0}^p(-1)^jr^{p-i}(1-r)^i\delta_{i,0}\zeta(\{2\}^j)\zeta^\star(\{2\}^{p-i-j})=T\sum\limits_{j=0}^p(-1)^jr^{p}\zeta(\{2\}^j)\zeta^\star(\{2\}^{p-j})=0.
\end{align*}
Hence, we have
\begin{align*}
  &\zeta^{T,r}_\ast(\{2\}^p,1)=2\sum\limits_{n=1}^p(-1)^n\zeta(2n+1)D_1(n,p)+2\sum\limits_{n=1}^p(-1)^{n+1}(1-2^{-2n})\zeta(2n+1)D_2(n,p),
\end{align*}
where
\begin{align*}
  D_1(n,p)&=\sum\limits_{i+j=0\atop i,j\geq0}^p(-1)^jr^{p-i}(1-r)^i\binom{2n}{2i}\zeta(\{2\}^{i+j-n})\zeta^\star(\{2\}^{p-i-j})\zeta(2n+1),\\
  D_2(n,p)&=\sum\limits_{i+j=0\atop i,j\geq0}^p(-1)^jr^{p-i}(1-r)^i\binom{2n}{2j+1}\zeta(\{2\}^{i+j-n})\zeta^\star(\{2\}^{p-i-j})\zeta(2n+1).
\end{align*}
For $D_1(n,p)$, we get
\begin{align*}
  D_1(n,p)&=\sum\limits_{i=0}^n\sum\limits_{j=0}^{p-n}(-1)^{n-i+j}r^{p-i}(1-r)^{i}\binom{2n}{2i}\zeta(\{2\}^j)\zeta^\star(\{2\}^{p-n-j})\\
  &=\sum\limits_{i=0}^n(-1)^{n-i}r^{p-i}(1-r)^i\binom{2n}{2i}\sum\limits_{j=0}^{p-n}(-1)^j\zeta(\{2\}^j)\zeta^\star(\{2\}^{p-n-j}).
\end{align*}
Using \eqref{Eq:zeta and t k} and \eqref{Eq:Gouldbinomal2}, we have
\begin{align*}
  D_1(n,p)=\begin{cases}
    (-1)^p\Omega_r^p & n=p,\\
    0 & 1\leq n\leq p-1.
  \end{cases}
\end{align*}
Similarly, for $D_2(n,p)$, we get
\begin{align*}
  D_2(n,p)&=\sum\limits_{j=0}^{n-1}(-1)^jr^j(1-r)^{n-j}\binom{2n}{2j+1}\sum\limits_{i=0}^{p-n}r^{p-n-i}(1-r)^i
\zeta(\{2\}^i)\zeta^\star(\{2\}^{p-n-i}).
\end{align*}
By \eqref{Eq:zeta^r k} and \eqref{Eq:Gouldbinomal1}, we have
\begin{align*}
  D_2(n,p)=(1-r)(-1)^{n-1}\Theta_r^n\zeta^r(\{2\}^{p-n}).
\end{align*}
Hence, we obtain
\begin{align*}
  \zeta^{T,r}_\ast(\{2\}^p,1)=2\sum\limits_{n=1}^p(1-r)(1-2^{-2n})\Theta_r^n\zeta^r(\{2\}^{p-n})\zeta(2n+1)+2\Omega_r^p\zeta^T_\ast(2p+1).
\end{align*}
Therefore, we complete the proof.
\qed

Taking $r=0$, \eqref{Eq:rmzvs-two-one} and \eqref{Eq:rmzvs-one-two} are special cases of \eqref{Eq:rmzvs-two-one-two} when $q=0$ and $p=0$. Taking $r=1$, \eqref{Eq:rmzvs-two-one} shows that $\zeta^\star(\{2\}^p,1)=2\zeta(2p+1)$ for $p\geq1$ and \eqref{Eq:rmzvs-one-two} shows that for $q\geq0$,
\begin{align*}
\zeta^{T,\star}_\ast(1,\{2\}^q)=2\sum\limits_{n=1}^{q}\zeta^\star(\{2\}^{q-n})\zeta(2n+1)-4q(1-2^{-2q})\zeta^T_\ast(2q+1)+T\zeta^\star(\{2\}^q).
\end{align*}
Moreover, we have the following evaluation formulas when $r=\frac{1}{2}$. 

\begin{cor}
  For any nonnegative integers $p$ and $q$, we have
  \begin{align*}
    \zeta^{T,\frac{1}{2}}_\ast(\{2\}^p,1)&=2\sum\limits_{n=1}^{\lfloor (p+1)/2\rfloor}(-1)^{n-1}(1-2^{-4n+2})\zeta^{\frac{1}{2}}(\{2\}^{p-2n+1})\zeta(4n-1)\\&\quad+2\cos\frac{p\pi}{2}\zeta^T_\ast(2p+1)-\delta_{p,0}T,\\
    \zeta_{\ast}^{T,\frac{1}{2}}(1,\{2\}^q)&=2\sum\limits_{n=1}^{\lfloor q/2\rfloor}(-1)^n\zeta^{\frac{1}{2}}(\{2\}^{q-2n})\zeta(4n+1)\\&\quad-2(1-2^{-2q})\sin\frac{q\pi}{2}\zeta^T_\ast(2q+1)+T\zeta^{\frac{1}{2}}(\{2\}^q),
  \end{align*}
  where $\lfloor k\rfloor$ is the largest integer not greater than $k$.
\end{cor}
\proof
Notice that 
\begin{align}\label{Eq:Theta_half}
  \Theta_{1/2}^n=2\sin\frac{n\pi}{2}=\begin{cases}
    0 & n:\text{even},\\
    2(-1)^{\frac{n-1}{2}} & n:\text{odd},
  \end{cases}
\end{align}
and
\begin{align}\label{Eq:Omega_half}
  \Omega_{1/2}^n=\cos\frac{n\pi}{2}=\begin{cases}
    (-1)^{\frac{n}{2}} & n:\text{even},\\
    0 & n:\text{odd}.
  \end{cases}
\end{align}
Hence we get the desired results.
\qed

Now we deal with the interpolated multiple $t$-values.
\begin{thm}\label{Thm:rmtvs-two-three}
For any nonnegative integers $p$ and $q$, we have
\begin{align}
t^r(\{2\}^p,3)&=\sum\limits_{n=0}^{p}2^{-2n-2}(\delta_{n,p}(1-2^{-2p-2})+1)\Theta_r^{n+1}t^r(\{2\}^{p-n})\zeta(2n+3),\label{Eq:r-mtv-l}\\
t^r(3,\{2\}^q)&=\sum\limits_{n=0}^{q}2^{-2n-2}(\delta_{n,q}+1-2^{-2n-2})\Theta_r^{n+1}t^r(\{2\}^{q-n})\zeta(2n+3)\label{Eq:r-mtv-r}.
\end{align}
\end{thm}
\proof
We prove \eqref{Eq:r-mtv-r}. Since $t(n)=(1-2^{-n})\zeta(n)$, \eqref{Eq:r-mtv-r} holds for $q=0$. For $q\geq1$, applying the map $\mathrm{t}$ to \eqref{Eq:S^r-r} with $a_1=z_3$ and $a_2=z_2$, we have 
\begin{align*}
t^r(3,\{2\}^q)=\sum\limits_{i+j=0\atop i,j\geq0}^q(-1)^ir^{q-j}(1-r)^jt(\{2\}^i,3,\{2\}^j)t^\star(\{2\}^{q-i-j}).
\end{align*}
As shown in \cite[Theorem 3]{Murakami2021}
\begin{align}
&t(\{2\}^i,3,\{2\}^j)\notag\\&=\sum\limits_{n=1}^{i+j+1}(-1)^{n-1}2^{-2n}\left\{(1-2^{-2n})\binom{2n}{2i+1}+\binom{2n}{2j+1}\right\}t(\{2\}^{i+j-n+1})\zeta(2n+1),\label{Eq:mtv2-3-2}
\end{align}
then we have
\begin{align*}
&t^r(3,\{2\}^q)
=\sum\limits_{n=1}^{q+1}(-1)^n2^{-2n}(1-2^{-2n})\zeta(2n+1)A_1(n,q)+\sum\limits_{n=1}^{q+1}(-1)^n2^{-2n}\zeta(2n+1)A_2(n,q),
\end{align*}
where
\begin{align*}
A_1(n,q)&=\sum\limits_{i+j=n-1\atop i,j\geq0}^q(-1)^{i-1}r^{q-j}(1-r)^{j}\binom{2n}{2i+1}t(\{2\}^{i+j-n+1})t^\star(\{2\}^{q-i-j}),\\
A_2(n,q)&=\sum\limits_{i+j=n-1\atop i,j\geq0}^q(-1)^{i-1}r^{q-j}(1-r)^{j}\binom{2n}{2j+1}t(\{2\}^{i+j-n+1})t^\star(\{2\}^{q-i-j}).
\end{align*}
For $A_1(n,q)$, we have
\begin{align*}
A_1(n,q)&=\sum\limits_{i=0}^{n-1}\sum\limits_{j=n-1-i}^{q-i}(-1)^{i-1}r^{q-j}(1-r)^{j}\binom{2n}{2i+1}t(\{2\}^{i+j-n+1})t^\star(\{2\}^{q-i-j})\\
&=\sum\limits_{i=0}^{n-1}(-1)^{i-1}r^{i}(1-r)^{n-1-i}\binom{2n}{2i+1}\sum\limits_{j=0}^{q-n+1}r^{q-n+1-j}(1-r)^jt(\{2\}^j)t^\star(\{2\}^{q-n+1-j}).
\end{align*}
According to \eqref{Eq:t^r k} and the combinatorial identity \eqref{Eq:Gouldbinomal1}, we get
\begin{align*}
A_1(n,q)&=\sum\limits_{i=0}^{n-1}(-1)^{i-1}r^{i}(1-r)^{n-1-i}\binom{2n}{2i+1}t^r(\{2\}^{q-n+1})\\
&=(-1)^n\Theta_r^nt^r(\{2\}^{q-n+1}).
\end{align*}
Similarly, for $A_2(n,q)$, we have
\begin{align*}
A_2(n,q)&=\sum\limits_{j=0}^{n-1}\sum\limits_{i=n-1-j}^{q-j}(-1)^{i-1}r^{q-j}(1-r)^{j}\binom{2n}{2j+1}t(\{2\}^{i+j-n+1})t^\star(\{2\}^{q-i-j})\\
&=\sum\limits_{j=0}^{n-1}(-1)^{n-j}r^{q-j}(1-r)^j\binom{2n}{2j+1}\sum\limits_{i=0}^{q-n+1}(-1)^it(\{2\}^i)t^\star(\{2\}^{q-n+1-i}).
\end{align*}
Using \eqref{Eq:zeta and t k} and \eqref{Eq:Gouldbinomal1}, we get
\begin{align*}
A_2(n,q)=\begin{cases}
(-1)^{q+1}\Theta_r^{q+1} & n=q+1,\\
0 & 1\leq n\leq q.
\end{cases}
\end{align*}
Therefore, we find
\begin{align*}
&t^r(3,\{2\}^q)=2^{-2q-2}\Theta_r^{q+1}\zeta(2q+3)+\sum\limits_{n=0}^{q}2^{-2n-2}(1-2^{-2n-2})\Theta_r^{n+1}t^r(\{2\}^{q-n})\zeta(2n+3),
\end{align*}
which implies the desired result.
\qed

Taking $r=0$, \eqref{Eq:r-mtv-l} and \eqref{Eq:r-mtv-r} are special cases of \eqref{Eq:mtv2-3-2} when $j=0$ and $i=0$. Taking $r=1$, \eqref{Eq:r-mtv-l} and \eqref{Eq:r-mtv-r} are special cases of \cite[Theorem 4.1]{Li-Wang}. And we get the following evaluation formulas from Theorem \ref{Thm:rmtvs-two-three} together with \eqref{Eq:Theta_half}.
\begin{cor}
For any nonnegative integers $p$ and $q$, we have
\begin{align*}
t^\frac{1}{2}(\{2\}^p,3)&=\sum\limits_{n=0}^{\lfloor p/2\rfloor }(-1)^n2^{-4n-1}\left(\delta_{2n,p}(1-2^{-4n-2})+1\right)t^\frac{1}{2}(\{2\}^{p-2n})\zeta(4n+3),\\
t^\frac{1}{2}(3,\{2\}^q)&=\sum\limits_{n=0}^{\lfloor q/2\rfloor}(-1)^n2^{-4n-1}\left(\delta_{2n,q}+1-2^{-4n-2}\right)t^\frac{1}{2}(\{2\}^{q-2n})\zeta(4n+3).
\end{align*}
\end{cor}

The following theorem gives evaluations for the harmonic regularized interpolated multiple $t$-values.
\begin{thm}\label{Thm:rmtvs-two-one}
For any nonnegative integers $p$ and $q$, we have
\begin{align}
t^{V,r}_\ast(\{2\}^p,1)&=\sum\limits_{n=1}^{p}2^{-2n}\left(1-2^{-2n}+\delta_{n,p}\right)\Omega_r^{n}t^r(\{2\}^{p-n})\zeta(2n+1)\notag\\
&\quad +\delta_{p,0}(V-\log2)+t^r(\{2\}^p)\log2,\label{Eq:r-mtv-regularized-1}\\
t^{V,r}_\ast(1,\{2\}^q)&=\sum\limits_{n=1}^{q}2^{-2n}\left(1+\delta_{n,q}(1-2^{-2n})\right)\Omega_r^{n}t^r(\{2\}^{q-n})\zeta(2n+1)\notag\\
&\quad+(V-\log2)t^r(\{2\}^q)+\delta_{q,0}\log2.\label{Eq:r-mtv-regularized}
\end{align}
\end{thm}
\proof
We prove \eqref{Eq:r-mtv-regularized}. It is obvious that \eqref{Eq:r-mtv-regularized} holds when $q=0$. For $q\geq1$, by \eqref{Eq:S^r-r}, we have 
\begin{align*}
t^{V,r}_\ast(1,\{2\}^q)=\sum\limits_{i+j=0\atop i,j\geq0}^q(-1)^ir^{q-j}(1-r)^jt^{V}_\ast(\{2\}^i,1,\{2\}^j)t^\star(\{2\}^{q-i-j}).
\end{align*}
Recall the evaluation
\begin{align}
t^{V}_\ast(\{2\}^i,1,\{2\}^j)&=\sum\limits_{n=1}^{i+j}(-1)^{n}2^{-2n}\left\{\binom{2n}{2i}+(1-2^{-2n})\binom{2n}{2j}\right\}t(\{2\}^{i+j-n})\zeta(2n+1)\notag\\
&\quad+\delta_{i,0}(V-\log2)t(\{2\}^j)+\delta_{j,0}t(\{2\}^i)\log2,\label{Eq:mtv2-1-2}
\end{align}
which was proved in \cite[(10)]{Charlton}. Then, we have
\begin{align*}
t^{V,r}_\ast(1,\{2\}^q)&=(V-\log2)\sum\limits_{i+j=0\atop i,j\geq0}^q(-1)^{i}r^{q-j}(1-r)^j\delta_{i,0}t(\{2\}^j)t^\star(\{2\}^{q-i-j})\\
&\quad+\log2\sum\limits_{i+j=0\atop i,j\geq0}^q(-1)^{i}r^{q-j}(1-r)^j\delta_{j,0}t(\{2\}^i)t^\star(\{2\}^{q-i-j})\\
&\quad+\sum\limits_{i+j=1\atop i,j\geq0}^q\sum\limits_{n=1}^{i+j}(-1)^{n+i}2^{-2n}r^{q-j}(1-r)^j\left\{\binom{2n}{2i}+(1-2^{-2n})\binom{2n}{2j}\right\}\\
&\quad\quad\times t(\{2\}^{i+j-n})t^\star(\{2\}^{q-i-j})\zeta(2n+1).
\end{align*}
Changing the order of summation and using \eqref{Eq:zeta and t k} and \eqref{Eq:t^r k}, we get
\begin{align*}{}
t^{V,r}_\ast(1,\{2\}^q)&=\sum\limits_{n=1}^q(-1)^n2^{-2n}\zeta(2n+1)C_1(n,q)+\sum\limits_{n=1}^q(-1)^n2^{-2n}(1-2^{-2n})\zeta(2n+1)C_2(n,q)\\
&\quad+(V-\log2)t^r(\{2\}^q)+\delta_{q,0}\log2,
\end{align*}
where
\begin{align*}
C_1(n,q)&=\sum\limits_{i+j=n\atop i,j\geq0}^q(-1)^ir^{q-j}(1-r)^j\binom{2n}{2i}t(\{2\}^{i+j-n})t^\star(\{2\}^{q-i-j}),\\
C_2(n,q)&=\sum\limits_{i+j=n\atop i,j\geq0}^q(-1)^ir^{q-j}(1-r)^j\binom{2n}{2j}t(\{2\}^{i+j-n})t^\star(\{2\}^{q-i-j}).
\end{align*}
For $C_1(n,q)$, we have
\begin{align*}
C_1(n,q)&=\sum\limits_{i=0}^n\sum\limits_{j=n-i}^{q-i}(-1)^ir^{q-j}(1-r)^j\binom{2n}{2i}t(\{2\}^{i+j-n})t^\star(\{2\}^{q-i-j})\\
&=\sum\limits_{i=0}^n(-1)^ir^i(1-r)^{n-i}\binom{2n}{2i}\sum\limits_{j=0}^{q-n}r^{q-n-j}(1-r)^jt(\{2\}^j)t^\star(\{2\}^{q-n-j}).
\end{align*} 
Using \eqref{Eq:t^r k} and \eqref{Eq:Gouldbinomal2}, we get
\begin{align*}
C_1(n,q)=(-1)^n\Omega_r^nt^r(\{2\}^{q-n}).
\end{align*}
Similarly, for $C_2(n,q)$, we have
\begin{align*}
C_2(n,q)&=\sum\limits_{j=0}^n(-1)^{n-j}r^{q-j}(1-r)^j\binom{2n}{2j}\sum\limits_{i=0}^{q-n}(-1)^it(\{2\}^{i})t^\star(\{2\}^{q-n-i})\\
&=\begin{cases}
(-1)^q\Omega_r^{q} & n=q,\\
0 & 1\leq n\leq q-1.
\end{cases}
\end{align*}
Hence, we obtain 
\begin{align*}
t^{V,r}_\ast(1,\{2\}^q)&=\sum\limits_{n=1}^q2^{-2n}\Omega_r^nt^r(\{2\}^{q-n})\zeta(2n+1)+2^{-2q}(1-2^{-2q})\Omega_r^{q}\zeta(2q+1)\\
&\quad+(V-\log2)t^r(\{2\}^q)+\delta_{q,0}\log2,
\end{align*}
which is the desired result.
\qed
 
In the case of $r=0$, \eqref{Eq:r-mtv-regularized-1} and \eqref{Eq:r-mtv-regularized} are special cases of \eqref{Eq:mtv2-1-2} when $j=0$ and $i=0$. And in the case of  $r=1$, \eqref{Eq:r-mtv-regularized-1} and \eqref{Eq:r-mtv-regularized} are special cases of \cite[Theorem 4.2]{Li-Wang}. Furthermore, we get the following evaluation formulas from Theorem \ref{Thm:rmtvs-two-one} together with \eqref{Eq:Omega_half}.
\begin{cor}
For any nonnegative integers $p$ and $q$, we have
\begin{align*}
t^{V,\frac{1}{2}}_\ast(\{2\}^p,1)&=\sum\limits_{n=1}^{\lfloor p/2\rfloor}(-1)^n2^{-4n}\left(1-2^{-4n}+\delta_{2n,p}\right)t^\frac{1}{2}(\{2\}^{p-2n})\zeta(4n+1)\\
&\quad +\delta_{p,0}(V-\log2)+t^\frac{1}{2}(\{2\}^p)\log2,\\
t^{V,\frac{1}{2}}_\ast(1,\{2\}^q)&=\sum\limits_{n=1}^{\lfloor q/2\rfloor}(-1)^n2^{-4n}\left(1+\delta_{2n,q}(1-2^{-4n})\right)t^\frac{1}{2}(\{2\}^{q-2n})\zeta(4n+1)\\
&\quad+(V-\log2)t^\frac{1}{2}(\{2\}^q)+\delta_{q,0}\log2.
\end{align*}
\end{cor}

We may regard Theorem \ref{Thm:r-mzvs2-3}, Theorem \ref{Thm:r-mzvs2-1}, Theorem \ref{Thm:rmtvs-two-three} and Theorem \ref{Thm:rmtvs-two-one} as ``weak'' generalizations of evaluations of MZ(S)Vs and M$t$(S)Vs. The evaluations for $\zeta^r(\{2\}^p,3,\{2\}^q)$ and $\zeta^r(\{2\}^p,1,\{2\}^q)$ are still unknown when $p,q\geq1$, as well as for $r$-M$t$Vs. Considering $r = \frac{1}{2}$, we use the software PARI/GP and get some identities when $p,q \geq 1$ and $p+q \leq 3$.  Here are the obtained identities:
\begin{align*}
  &\zeta^\frac{1}{2}(2,3,2)=-5\zeta(2)\zeta(5)+\frac{73}{8}\zeta(7),\\
  &\zeta^\frac{1}{2}(2,3,2,2)=2\zeta^\frac{1}{2}(2,2)\zeta(5)+\frac{115}{8}\zeta(2)\zeta(7)-\frac{105}{4}\zeta(9),\\
  &\zeta^\frac{1}{2}(2,2,3,2)=-7\zeta^\frac{1}{2}(2,2)\zeta(5)-\frac{155}{16}\zeta(2)\zeta(7)+\frac{105}{4}\zeta(9),\\
  &\zeta^\frac{1}{2}(2,1,2)=-\frac{3}{2}\zeta(2)\zeta(3)+\frac{17}{4}\zeta(5),\\
  &\zeta^\frac{1}{2}(2,1,2,2)=\frac{3}{5}\zeta^\frac{1}{2}(2,2)\zeta(3)+\frac{11}{2}\zeta(2)\zeta(5)-\frac{301}{32}\zeta(7),\\
  &\zeta^\frac{1}{2}(2,2,1,2)=-\frac{21}{10}\zeta^\frac{1}{2}(2,2)\zeta(3)-\frac{13}{4}\zeta(2)\zeta(5)+\frac{301}{32}\zeta(7),\\
  &t^\frac{1}{2}(2,3,2)=\frac{381}{1024}\zeta(7),\\
  &t^\frac{1}{2}(2,3,2,2)=\frac{5}{48}t^\frac{1}{2}(2,2)\zeta(5)+\frac{87}{1024}t(2)\zeta(7),\\
  &t^\frac{1}{2}(2,2,3,2)=-\frac{5}{48}t^\frac{1}{2}(2,2)\zeta(5)+\frac{461}{2048}t(2)\zeta(7),\\
  &t^\frac{1}{2}(2,1,2)=\frac{31}{64}\zeta(5),\\
  &t^\frac{1}{2}(2,1,2,2)=\frac{1}{8}t^\frac{1}{2}(2,2)\zeta(3)+\frac{29}{256}t(2)\zeta(5),\\
  &t^\frac{1}{2}(2,2,1,2)=-\frac{1}{8}t^\frac{1}{2}(2,2)\zeta(3)+\frac{1}{4}t(2)\zeta(5).
\end{align*}

However, using PARI/GP, when $p,q\geq1$ and $p+q=4$, we have not found such linear relations among $\zeta^\frac{1}{2}(\{2\}^p,3,\{2\}^q)$ and $\zeta^\frac{1}{2}(\{2\}^{p+q-n+1})\zeta(2n+1)(1\leq n\leq 5)$. Also, there are no linear relations among $\zeta^\frac{1}{2}(\{2\}^p,1,\{2\}^q)$ and $\zeta^\frac{1}{2}(\{2\}^{p+q-n})\zeta(2n+1)(1\leq n\leq 4)$. Moreover, such explicit linear relations are also unavailable for the cases of $r$-M$t$Vs involving $t^\frac{1}{2}(\{2\}^p,3,\{2\}^q)$ and $t^\frac{1}{2}(\{2\}^p,1,\{2\}^q)$ when $p,q\geq1$ and $p+q=4$.

\section*{Acknowledgments}
The authors would like to thank Prof. Masanobu Kaneko for his tutorial on PARI/GP. This work was supported by the International Exchange Program for Graduate Students, Tongji University (grant number 2023020002).

\end{document}